\documentclass[11pt,twoside]{article}
\usepackage{amsfonts,amssymb}

\def\Bel{\hbox{\vrule\vbox to 7pt
{\hrule width 6pt
\vfill\hrule}\vrule}} 
\def\Bbb{\mathbb}
\def\qed{\space\Bel}
\def\goth{\mathfrak}
\def\cal{\mathcal}
\newtheorem{theorem}{Theorem}
\newtheorem{lemma}{Lemma}
\newtheorem{conjecture}{Conjecture}
\newtheorem{example}{Example}

\title{Note on Trace Class Groups}
\author{Gerrit van Dijk}
\date{}
\begin{document}
\maketitle

\begin{abstract} \noindent
A Lie group $G$ is called a trace class group if for every irreducible unitary representation $\pi$ of $G$ and every $C^\infty$ function $f$ with compact support the operator $\pi (f)$ is of trace class. In this note we prove that the semidirect product of $\Bbb R^n$and a real semisimple algebraic subgroup $G$ of ${\rm GL}(n,
\Bbb R)$ is a trace class group only if $G$ is compact. The converse has been shown elsewhere. We also make a descent start with the study of semidirect products with Heisenberg-type groups.
  \end{abstract}
\section{Introduction}
In this note we resume the study of trace class groups from [4]. An irreducible unitary representation $\pi$ of a Lie group is said to be of trace class if for every $C^\infty$ function $f$ with compact support the operator $\pi (f)$ is of trace class. A Lie group is said to be of trace class, or briefly, a trace class group, if every irreducible unitary representation is of trace class. Well-known examples of such groups are reductive Lie groups and unipotent Lie groups. In general each (real algebraic) Lie group is a semidirect product of a reductive and a unipotent Lie group. One of the highlights of [4] is the theorem that the semidirect product of a real algebraic semisimple Lie group and its Lie algebra is a trace class group if and only if the group is compact. In this note we prove a generalization of this theorem, provided by the case of a semisimple real algebraic group $G$ acting on a real finite-dimensional vector space $V$ by linear transformations and considering the semidirect product of $V$ and $G$. We also make a beginning with the study of real algebraic groups with unipotent radical equal to a Heisenberg group.
\section{Formula for the character of an induced representation}
Let $G$ be a locally compact group and $N$ a closed subgroup. Choose right Haar measures $dg$ on $G$ and $dn$ on $N$. We may find a strictly positive continuous function $q$ on $G$ satisfying
$$q(e)=1,$$
\begin{equation}q(ng)=\Delta_N(n)\,\Delta_G(n^{-1})\, q(g)\quad(n\in N, q\in G),\end{equation}
where $\Delta_N, \Delta_G$ denote the modular functions on $N, G$. For example
$$\int_G f(a^{-1}g) dg=\Delta_G (a)\,\int_G f(g) dg$$
for all $a\in G$ and $f\in C_c(G)$.\
\par\noindent
The function $q$ defines a quasi-invariant measure $d_q\dot{g}$ on $\dot{G}=G/N$ (the space of right cosets with respect to $N$) as follows. For $f\in C_c(G)$ set $T_Nf({\dot g})=\int_N f(ng)dn,\, {\dot g}=Ng$. Then $d_q({\dot g})$ is defined by
$$\int_{G/N} T_Nf({\dot g})d_q({\dot g})=\int_G f(g)q(g)dg.$$
Let $\gamma$ be a unitary representation of $N$ and set $\pi={\rm Ind}_N^G\gamma$. We write down a formula for the character of $\pi$ in terms of that for $\gamma$. Let us give the definition of $\pi$. Let $\cal H$ be the Hilbert space of $\gamma$. Then $\pi$ acts on the space ${\cal H}_\pi$ of function $f: G\to {\cal H}$ satisfying
$$f(ng)=\gamma (n) f(g)\quad {\rm and}\quad \int_{
G/N}\Vert f({\dot g})\Vert^2 d_q{\dot g}< \infty.$$
The action of $\pi$ is
$$\pi (g)f(x)=f(xg)[q(xg)/q(x)]^{1/2}.$$
\begin{theorem}{\rm ([2], Theorem 3.2)}.
Let $\varphi\in C_c(G), \varphi^\ast(g)=\overline{\varphi (g^{-1})}\Delta_G (g^{-1})$ and set $\psi=\varphi \ast_G \varphi^\ast$. Then
$${\rm tr}\,\pi (\psi)=$$
\begin{equation}\int_{G/N}\Delta_G(g^{-1})q(g^{-1})\,{\rm tr} [\int_N \psi (g^{-1}ng)\gamma (n)\Delta_G(n)^{1/2}\Delta_N(n)^{-1/2}dn ]\, d_q{\dot g}\end{equation}
in the sense that both sides are finite and equal or both $+\infty$.
\end{theorem} 
\noindent
A group $G$ is called unimodular if $\Delta_G =1$. If $G$ is a unimodualr Lie group we have ${\rm tr}\,\pi(\varphi)$ is finite for all functions $\varphi\in C_c^\infty(G)$ if and if ${\rm tr}\,\pi(\psi)$ is for all fuctions $\psi$ of the form $\psi=\varphi\ast_G \varphi^\ast$ with $\varphi\in C_c^\infty (G)$. This is because any $\varphi$ is a finite sum of functions of the form $\psi$ by [1].
\section{Application to semidirect products}
\setcounter{theorem}{0}
Let $V$ be a finite-dimensional real vector space and $H$ a closed subgroup of ${\rm GL}(V)$. Set $G=V\rtimes H$. The product in $G$ is given by 
$$(v,h)(v',h')=(v + h\cdot v', hh')\quad {\rm for}\  v,v'\in V, h,h'\in H.$$
If $dh$ is a right Haar measure on $H$ and $dv$ one on $V$, then $dg=dvdh$ is a right haar measure on $G$. Denote by ${\widehat V}$ the space of continuous unitary characters of $V$. For each $\chi\in\widehat V$ and each $h\in H$ consider the function 
$$v\mapsto \chi (h\cdot v)\quad (v\in V).$$
This is again a continuous unitary character of $V$, which we call $\chi\cdot h$. The set of all $\chi\cdot h$ is called the orbit of $\chi$ in $\widehat V$ and 
$$H_\chi=\{h\in H:\, \chi\cdot h=\chi\}$$
the stability subgroup of $\chi$. Choose an irreducible unitary representation $\rho$ of $H_\chi$ and define $\chi\otimes\rho$ by
$$(v,h)\mapsto \chi (v)\,\rho (h)\quad (v\in V,\, h\in H_\chi).$$
Then $\chi\otimes\rho$ is an irreducible unitary representation of $V\rtimes H_\chi$ and the induced representation $\pi_{\chi,\rho}$ is an irreducible unitary representation of $V\rtimes H$, see [3], p. 43. Now apply (2) with $G=V\rtimes H$, $N=V\rtimes H_\chi$ and $\pi=\pi_{\chi,\rho}$. Choose as before $dvdh$ and similarly $dvdh_\chi$ as right Haar mesures on $G$ and $N$. Then $\Delta_G$ and $\Delta_N$ are given by
$$\Delta_G(v,h)={\rm det}(h).\,\Delta_H(h)\quad (v\in V, h\in H)$$ and similarly 
$$\Delta_N(v,h)={\rm det}(h)\,\Delta_{H_\chi}(h)\quad (v\in V, h\in H_\chi ).$$
Notice that $q$ is left $V$-invariant, so we may write $q(v,h)=Q(h)$. $Q$ satisfies 
$$Q(h_0h)=\Delta_{H_\chi}(h_0)\Delta_H(h_0^{-1})Q(h)$$
for $h\in H,\, h_0\in H_\chi$. Let us now rewrite (2) for the above particular case.
\begin{theorem}
Let $\psi$ be as in Theorem 2.1. Then 
$${\rm tr}\, \pi(\psi)=\int_{H/{H_\chi}} \Delta_H(h)^{-1}\, Q(h)^{-1}.$$
\begin{equation}
{\rm tr}\, [\int_V\int_{H_\chi}\psi(v, h^{-1}h_oh)\, (\chi\cdot h)(v)\, \rho(h_0)\,\Delta_H(h_0)^{1/2}\,\Delta_{H_\chi}(h_0)^{-1/2}\, dvdh_0]\, d_Q{\dot h}.
\end{equation}
\end{theorem}
\section{A special case}
\setcounter{theorem}{0}
Let $V$ be a finite-dimensional real vector space with inner product $\langle\ ,\ \rangle$ and with complexification $\bf V$.  Denote by $\bf G$ a connected , complex, semisimple, linear algebraic subgroup of ${\rm GL}(\bf V)$. Assume that $\bf G$ is defined over $\Bbb R$ and set $G={\bf G}(\Bbb R)$ for its  group of real points. Then $G$ is a semisimple Lie group with finite center and finitely many connected components. For any $g\in {\rm GL}(V)$ set $\langle g\cdot v,\, w\rangle=\langle v,\, ^t g\cdot w\rangle$ for all $v,w\in V$. Assume that $G$ is invariant under the Cartan involution $\theta$ defined by $\theta (g)=\,^tg^{-1}\ (g\in G)$. Write $G_1=V\rtimes G$. Notice that $G_1$ is unimodular. The purpose of this note is to show 
\begin{theorem} The group $G_1$ is a trace class group only if the group $G$ is compact.
\end{theorem}
The converse of this theorem has been proved in [4], Lemma 14.2.
\vskip.2cm\noindent
{\bf Proof}.
Assume $G$ to be non-compact. The Cartan involution of $G$ gives rise to a Cartan involution of the Lie algebra $\goth g$ of $G$, that we again denote by $\theta$. Write $\goth g=\goth k +\goth p$ for the decomposition of $\goth g$ into $\pm 1$-eigenspaces of $\theta$. Then $\goth k$ consists of anti-symmetric and $\goth p$ of symmetric elements. Set $K=G\cap {\rm O}(n,\Bbb R)$.  Then $\goth k$ is the Lie algebra of $K$. Select a non-trivial maximal Abelian subspace $\goth a$ of $\goth p$, which exists because $G$ is non-compact, and let $\Sigma$ denote the set of roots of $(\goth g,\,\goth a)$. Then $\Sigma$ is root system (with multiplicities). Let $\Delta$ be a set of simple roots and $\Sigma^+$ the set of positive roots with respect to $\Delta$. Denote by $\goth n$ the Lie subalgebra spanned by the positive root vectors and by $N$ the corresponding algebraic subgroup. Then one has $\goth g=\goth k +\goth a +\goth n$ and similarly $G=KAN$, the Iwasawa decomposition of $G$. Let $\xi_0$ be a highest weight vector in $V$ with highest weight $\lambda\not= 1$ (with respect to $N$ and $A$). Such a vector exists. Indeed, if all highest weight vectors have weight is equal to one, then $G=\{1\}$, which is not the case. Set $G_0$ for the stabilizer of $\xi_0$ in $G$. Denote by $B$ the Killing form of $\goth g$ and define $H_1$ by
$$\lambda (H)=B(H,\, H_1)\quad (H\in\goth a),$$ 
and set $A_1={\rm exp}\, \Bbb R H_1$.  Let us denote by $P$ the stabilizer of the half-line ${\Bbb R}^\ast_+\,\xi_0$. Then $P=A_1G_0$, where $A_1\cap G_0=\{1\}, \, a_1G_0a_1^{-1}=G_0$ for all $a_1\in A_1$. Let $da_1$ and $dg_0$ denote right Haar measures on $A_1$ and $G_0$ respectively. Then $dp=\delta (a_1)da_1dg_0$, where $\delta(a_1)={\rm det}\, {\rm Ad}(a_1^{-1})|_{\goth g_0}$ is a right Haar measure  on $P$. Since $A\subset P$ and $N\subset P$ we have $G=PK$ and
$$dg=d(g_0a_1k)=\Delta_{G_0}(g_0^{-1})\,\delta(a_1^{-1})\, dg_0da_1dk.$$
Define the character $\chi_0$ of $V$ by
$$\chi_0(v)={\rm e}^{-2\pi i\langle \xi_0,\,,v\rangle}\quad (v\in V).$$
Then ${\rm Stab}\,\chi_0=\theta (G_0)$. Let us write $H_0=\theta (G_0)$, and let $dh_0$ be a right Haar measure on $H_0$. In a similar way as above we have
$$dg=d(h_0a_1k)=\Delta_{H_0}(h_0)\delta_0(a_1)\,dh_0da_1dk,$$
where $\delta_0(a_1)={\rm det}\,{\rm Ad}(a_1)|_{\goth h_0}.$
We will consider the representation $\pi$ given by
$$\pi={\rm Ind}\,_{V\rtimes H_0}^{V\rtimes G}\,\chi_0\otimes  1$$
and determine whether its trace exists. Let
$$Q(h_0a_1k)=\Delta_{H_0}(h_0)\, \delta_0(a_1).$$
Equation (3) then becomes:
\begin{equation} 
{\rm tr}\, (\psi)=\int_{A_1}\widehat\psi_1(a_1\cdot \xi_0)\,[\, \int_{H_0}\psi_2(h_0)\,\Delta_{H_0}(h_0)^{-1/2}\, dh_0\,]\, da_1
\end{equation}
where we toke $\psi\in C_c^\infty (G)$ of the form $\psi (v,g)=\psi_1(v)\,\psi_2(g)\ (v\in V,g\in G)$, both $\psi_1$ and $\psi_2$ $K$-invariant and $\widehat\psi_1(w)=\int_V \psi_1(v)\,{\rm e}^{-2\pi i \langle v,\, w\rangle} dv$ for $w\in V$. Clearly the integral 
$$\int_{A_1}\widehat\psi_1(a_1\cdot \xi_0)\, da_1=\int_0^\infty \widehat\psi_1(\mu\xi_0)\, {d\mu\over \mu}$$ 
diverges for suitable $\psi_1$. So we may conclude that $G$ cannot be non-compact. This concludes the proof of the theorem.
\hfill\qed
\section{Semidirect products with Heisenberg groups}
In this section we extend our scope to semidirect products $G=V\rtimes H$ with $V$ a non-necessarily Abelian normal subgroup of $G$. Notice that any real algebraic group is of this form according to the Levi decomposition, see [4], Proposition 2.1.
Let us begin with some preparations.
\begin{lemma} Let $G$ be a Lie group and $N$ a closed normal subgroup of $G$. If $G$ is a trace class group, then the quotient group $\dot{G}=G/H$ is
\end{lemma} 
{\bf Proof}. Let $\dot{\pi}$  be an irreducible unitary represenatation of $\dot{G}$ and $\pi$ the corresponding representation of $G$. Choose a right Haarmeasure $dn$ on $N$ and define for $f\in C_c^\infty (G)$
$$T_Nf(\dot{x})=\int_{G/N} f(nx)dn.$$
Then $f\mapsto T_Nf$ is a continuous surjective linear map $C_c^\infty(G)\to C_c^\infty (\dot{G})$ and one has
$$\pi (f)=\dot{\pi}(T_Nf)\quad (f\in C_c^\infty (G)).$$
So the result follows.
\hfill\qed
\vskip.1cm\noindent
Let us consider a special case. Denote by $G$ a Lie group, being the semidirect product $G=V\rtimes H$ where $H$ is a closed subgroup of $G$ and $V$ a closed normal subgroup of $G$. Let $W$ be a closed normal subgroup of $G$ contained in $V$. Denote by $v\mapsto\dot{v}$ the canonical map $V\to V/W=\dot V$. The group $H$ acts on $\dot{V}$ by
$$h\cdot \dot{v}=\dot{\widetilde{h\cdot v}}\quad (h\in H, v\in V).$$
Set $\dot{G}=\dot{V}\rtimes H$. Then we have
\setcounter{theorem}{2}
\begin{lemma} The map $(v,h)\mapsto (\dot{v}, h)$ is a surjective homomorphism from $G$ to $\dot{G}$ with kernel $W$.
\end{lemma}
\par\noindent
We will now specialize to real algebraic groups $G$ of the form $G=V\rtimes H$ with $V$ a unipotent and $H$ a semisimple real algebraic group. Set
$$V_1\supset V_2\supset \cdots \supset V_k\supset \{1\}$$
for the descending series of $V$, where $V_1=V,\, V_i=[V,V_{i-1}],\, V_k=Z$, the (non-trivial center of $V$. Notice that each $V_i$ is normal in $G$, so in particular $H$-invariant. Define 
$$L=\{h\in H:\, h=id\ {\rm on}\ V_i/V_{i+1}\ {\rm for\ all}\ i=1,\ldots, k\}$$
Clearly $L$ is a closed real algebraic normal subgroup of $H$. We can now formulate a conjecture.
\setcounter{theorem}{3}
\begin{conjecture}
Let $G=V\rtimes H$ with $V$ a unipotent and $H$ a semisimple real algebraic group. Then $G$ is a trace class group if and only if $H/L$ is compact.
\end{conjecture}
Let us consider an example with $V$ of a special nature, namely a Heisenberg group. Such a group can be seen as the most simple choice for a non-Abelian group $V$. We shall show that the conjecture holds in this case.
\setcounter{theorem}{4}
\begin{example}\rm
Denote by $\cal V$ the $(2n+1)$-dimensional Heisenberg group with Lie algebra basis $x_1,\ldots ,x_n, y_1,\ldots ,y_n, z$ with $[x_i,\, y_i]=z$ and all other brackets equal to zero. Then $Z$ is one-dimensional and spanned by $z$. Let $H$ be any semisimple real algebraic group acting on $\cal V$ algebraically and set $G={\cal V}\rtimes H$. There are two kinds of irreducible unitary representations $\pi$ of $\cal V$, depending on their behaviour on $Z$. By Schur's Lemma we have $\pi (z)=\lambda (z)\, I\ (z\in Z)$ for some character $\lambda$ of $Z$.
\par\noindent
If $\lambda\not= 1$, then $\pi$ is equivalent with an infinite-dimensional representation $\pi_\lambda$ satisfying $\pi_\lambda (z)=\lambda (z)\, I$ and $\pi_\lambda$ is square-integrable modulo $Z$.
\par\noindent
If $\lambda =!$, then $\pi$ is actually a one-dimensional representation of ${\cal V}/Z\simeq {\Bbb R}^{2n}$, so a character $\chi$ of $\cal V$.
\par\noindent
Let us now perform the usual construction for the determination of the irreducible unitary representations of $G$, see [2], p. 470. Let us begin with $\pi_\lambda$. Since $H$ acts trivially on $Z$, one has ${\rm Stab}\, \pi_\lambda =H$. If $\rho$ is an irreducible unitary representation of $H$ then $\pi_\lambda\otimes \rho$ in one of $G$, of trace class. Let now $\chi$ be a character of ${\cal V}/Z$ (hence of $\cal V$). By the ususal construction (see [4]), we always obtain a trace class representation of ${\cal V}/Z\rtimes H=G/Z$, so of $G$, if and only if $H/L_0$ is compact, where $L_0=\{h\in H:\, h=id\ {\rm on}\ {\cal V}/Z\}$. Clearly $L_0=L$, so if and only if $H/L$ is compact. Resuming, $G={\cal V}\rtimes H$ is trace class if and only if $H/L$ is compact.
\end{example}
   
\vskip.2cm
{\small Gerrit van Dijk
\par\noindent
Mathematisch Instituut
\par\noindent
Niels Bohrweg 1
\par\noindent
2333 CA Leiden, The Netherlands
\par\noindent
dijk@math.leidenuniv.nl
\end{document}